\RequirePackage{fix-cm}

\documentclass[smallextended]{svjour3}       

\smartqed  

\usepackage{graphicx}
\usepackage{epstopdf}

\usepackage{amsmath}
\usepackage{mathrsfs}

\usepackage{algorithm}
\usepackage{algorithmic}

\usepackage{booktabs}
\usepackage{enumerate}
\usepackage{amsfonts}

\usepackage{url} 
\usepackage{lineno}
\usepackage{soul}

\begin{document}


\title{An integer programming model for the selection of pumped-hydro storage projects}
\author{Tiago Andrade \and
        Rafael Kelman \and
        Tain{\'a} M. Cunha \and
        Luiz Rodolpho de Albuquerque \and
        Rodrigo F. Calili}


\institute{T. Andrade, R. Kelman, T. M. Cunha, L. R. Albuquerque are with PSR at:\at
              Praia de Botafogo, 370 - Botafogo, Rio de Janeiro - RJ, 22250-040, Brazil \\
              Tel.: +55-21-3906-2100
              \\
              R. f. Calili is with PUC-Rio at:\at
              Marquês de São Vicente, 225 - Gávea, Rio de Janeiro - RJ, 22451-900, Brazil \\
              \hfill
              The authors would like to thank EDF, CTG, Brookfield and Light for jointly supporting and funding this research as part of a R\&D project. The project is supervised by ANEEL - the Brazilian electricity regulator.
              \\
                 \hfill
                The authors would also like to thank former consultant and chairman of the board of The Nature Conservancy, Mr. David Harrison, for his insights and comments.  
                \\
                \hfill
                The data used for case study in the study are available at Zenodo repository via  \url{https://doi.org/46510.5281/zenodo.4019256}
                with  Creative Commons Attribution 4.0 International \cite{andrade_tiago_2020_4019256}.
            \and
              Corresponding author e-mail: tiago.andrade@psr-inc.com
}

\date{April 20, 2020}

\maketitle

\begin{abstract}
Energy storage systems - in particular, Pumped Hydropower Storage (PHS) - will be increasingly important to support the transition of power systems toward zero emissions. The reason is that PHS can mitigate the variability and uncertainty of renewable energy production from solar and wind power to balance electricity demand with supply. In this paper, we propose an integer programming problem for PHS siting that uses a Digital Elevation Model (DEM) to meet an energy storage requirement. It assumes the existence of a reservoir, lake, or river, and decides where to build a reservoir that will constitute the PHS with the existing body of water. This model finds minimum-cost project candidates given parameters such as desired head, power, and operation time. The paper discusses different solution methods to assure reservoir closure and avoid its fragmentation. A heuristic explores the representations of the DEM, from more aggregate to more precise, to sequentially refine the solution based on the last selected site, which reduces computational effort. The formulation is general and the objective function includes both construction and equipment costs. Constraints are related to the energy storage target and reservoir closure based on the DEM. We illustrate the methodology by selecting multiple PHS projects next to the reservoir of the Sobradinho hydropower plant in Brazil. The result of this model can be seen as a bottom-up step that prepares PHS candidate projects to be considered by an integrated resource planning model in a top-down step, that would select from these shortlisted projects.
\end{abstract}



\section{Introduction}

At this time the world faces the dual challenge of increasing the availability of electrical energy while reaching toward a zero-carbon future in a way that minimizes any negative social and environmental impacts. The debate in the electrical power sector has taken a sharp turn recently. It is now about how fast carbon emissions can be reduced entirely, and how to do it. It is no longer relevant whether it is realistic, or even necessary, to slow the emission of greenhouse gases, and whether much of the world outside the developed countries could afford it. The conversation now is about the energy transition, from a hydrocarbon economy to a no-carbon one \cite{kelman:hal-02147740}.

The contribution of Variable Renewable Electricity (VRE), such as wind and solar photovoltaic projects, is increasing worldwide due to favorable conditions such as economic performance and resource availability \cite{BARBOUR2016421,CONNOLLY2016}. Technically gas fired plants could provide capacity to supply peak demand and reserves, but at the expense of increasing the emission of greenhouse gases. While conventional hydropower can also support the growth of VRE in power systems, concerns related to their socio-environmental impacts often make this alternative less appealing. The most competitive alternative is another hydraulic flexible resource: Pumped Hydropower Storage (PHS)  \cite{GUITTET2016560,IHA2018}. From 1960 to 1990, PHS development was mainly used to provide a flexible generation in systems with large share of base loaded nuclear or coal power in Europe, United States and Japan. More recently, the renewed interest on PHS is as a flexible renewable resource with several possible roles in more diversified electricity matrices, mostly in the USA, European countries and China \cite{BARBOUR2016421,ZHANG2015}.   

Energy storage can support energy security and climate change challenges by providing valuable benefits in power systems \cite{book}. As the cost of emerging technologies decreases, storage will become even more competitive. Meanwhile, the range of economical services granted by such solution will expand \cite{IRENA2017}. Compared to other energy storage technologies, PHS offers a longer service life, comparable efficiency, good response time, an ability to follow load changes, maintain frequency and support voltage stability \cite{REHMAN2015586}. Its operation is characterized by the pumping of water from a lower reservoir or river for storage into an upper reservoir during low demand period to  produce electricity during high demand period. These systems can operate under different topologies, such as connecting existing reservoirs or construction of two new closed-loop reservoirs (e.g. with no water influxes in either end)   \cite{osti_1165460}. The principal issue, however, is that PHS is highly site-specific and feasibility depends greatly on topographic and geologic variables of each potential project. 

Most of the published works dealing with the search for pumped storage sites are based on the guidelines of a workshop organized by the Institute for Energy and Transport (IET) of the Joint Research Center (JRC). The discussions were consolidated in a technical report \cite{SETIS}, that proposes general criteria considering the use of GIS-based tools and variables such as: the dimensions of existing or new reservoirs, the distance between potential reservoirs or other bodies of water, the head provided by the difference in water level between the two identified locations, the average slope of the land surface, the proximity to the electrical and the road network, in addition to socio-environmental impacts.

Studies developed after this workshop examine specific topologies (Table \ref{tab:topology}), adapting the reference values indicated by the IET to the characteristics of the chosen region and the extent of the assessment in spatial terms (local, regional, national or integrated system). These peculiarities define the scope of the search for potential locations in different studies, generally, limiting it to a given plant size, associated with the dimensions of the upper reservoir and a specific operating cycle.

\begin{table}[htbp]
  \footnotesize
  \centering
  \caption{PHS Topologies}
    \begin{tabular}{c l}
    \toprule
Topology & Description \\
\cmidrule{2-2}
T1 & Two existing reservoirs linked with one or several penstocks with a \\ & powerhouse to transform them to a PHS scheme \\
T2 & One existing lake or reservoir linked to a new reservoir built on a flat area, \\ & with excavation and embankments, on a depression or in a valley \\
T3 & A greenfield PHS with new reservoirs built in valleys, closed by dams, \\ & depressions or hill tops \\
T4 & Sea based PHS, in which the sea can be used as the lower reservoir \\ & combined with a new one or an upper basin with a cavern as a lower reservoir \\
T5 & Multi-reservoir systems combining PHS and conventional plants \\
T6 & Scheme with a large river that can provide enough inflow for the PHS system \\ &  used as the lower reservoir \\
T7 & Use of an abandoned mine pit as a reservoir as in a T2 scheme \\
    \bottomrule
    \end{tabular}%
    \label{tab:topology}
\end{table}%

Another usual classification comes from a water resources management perspective. Closed-loop systems, where neither reservoir is part of a river course, correspond to the T3 type. In semi-open systems, one of the reservoirs (frequently the lower one) is part of a water body, such as in T2, T4 and T6. Open systems have both reservoirs as part of a river course.

In 2012, a study in the Turkish territory \cite{FITZGERALD2012483} was based on an algorithm that seeks an upper reservoir to maximize the head between the two sites in a circular search region of 5 km surrounding an existing reservoir (T2 type). It also limited the search of the second reservoir to flat surfaces, considering that strongly sloped areas increase the complexity of the construction. A minimum reservoir area of 70 ha was considered, with 20 ha corresponding to civil structures and other components. Candidate sites were then ranked by energy storage capacity, in GWh, and the site of highest capacity related to an existing reservoir was selected.

SINTEF Energy Research developed a tool in 2013 \cite{ZINKE2013} to identify pumped-storage plants in Norway considering a pair of existing reservoirs (T1 type) and using a combination GIS resources and Python language. It defined a methodology built in 3 stages: topographic analysis, calculation of selection methods and choice of locations. The first step calculates the nearest distance between two existing sites and applies some geographic restrictions that are complemented by infrastructure and environmental proximity analysis in the second stage. For the selection of sites, three screening modes were defined, each one with a different parameter (power production, in MW, storage duration, in days, and water level change rate, in m/hour).

Recent studies explored other possibilities using GIS tools. In 2017, an investigation of pumped storage sites in Tibet \cite{LU20171045} searched for topology types T1 and T2 in a region with transmission network distant from potential sites. In the case of T2, the methodology takes advantage of a Multiple Ring Buffer tool around the existing site to intersect drainage lines while assuming a pyramidal geometry for the upper reservoir. 

A French screening \cite{ROGEAU2017241}, also carried out in 2017, added to these two types of scheme the closed cycle one (T3) to find small plants. This study establishes a method equivalent to the one applied in Norway. This paper considers that the small-PHS scheme can involve only two reservoirs that cannot be linked more than once. The strategy aimed at the definition of the maximum water flow for each connection and selects the pairs of reservoirs with higher head, as it can provide more installed capacity.

In 2018, researchers carried out a broad study for Iran \cite{GHORBANI2019854}, covering the topologies T1, T2, T4 and T6 using GIS tools combined with multicriteria analysis. Using spatial raster tools, the study considers the head ($H$), conveyance length-head ratio ($L/H$) and slope of ground surface, beyond other parameters, such as geological conditions, to select the best sites.

An Australian research \cite{LU2018300} conducted in the same year a worldwide survey to identify type T3 plants (closed-loop system), involving advanced geoprocessing algorithms considering two types of upper reservoir: dry-gully (DG), which seeks to close a topographic saddle with a dam in a mountainous topography, and turkey's nest (TN), which are equivalent to huge artificial pools, requiring the embankments in all its contour. This study resulted in an interesting search tool \cite{AREMI}, available for free on the internet. For DG schemes, it calculates the storage capacity of the reservoirs and the required earthwork to build a dam. It uses a raster watershed model from which a flow direction raster is derived. Within this watershed, dam wall heights are created to define the reservoirs, avoiding steep terrains by the exclusion of pour points located on a slope greater than 1:5. The dam contour is defined from the common cells of the watershed and its corresponding flooded areas. In TN schemes, it assumes that the soils or rocks required to build the dam are obtained from excavation within the reservoir area and hence adds half of the dam volume to the total volume of reservoir. It uses also a simplified form of dam volume calculation, in which freeboard and dam crest width are not considered. 

A GIS software named STORES was then developed using Python scripts comprising functional modules to exclude the regions without enough height within an acceptable distance (minimum $H/L$ ratio of 1:10), identify DG and TN sites within the highlighted areas and optimize the site selection. In this study, a ratio between the water storage capacity ($W$) divided by the earthwork required to build such a reservoir/ dam capable of this storage capacity ($R$) was used to select optimal sites, while cost models are not integrated in it.

Lately, a group of researchers developed a methodology to estimate the global pumped storage potential considering T6 schemes (open system) for a reservoir operation to regulate river flows  through a reduction of seasonality and inter-annual variations \cite{HUNT2020}. The screening begins by examining rivers with reasonable flow rate, being no more than 30 km away from an existing reservoir. Then, it verifies if a dam can be built around this area in a so-called Point Under Analysis (PUA), testing four cardinal orientations. Usually located in a tributary river, from this PUA, the model considers five heights (from 50m to 250m, at 50 m intervals), before removing projects with competing dams. Finding the flooded side of the dam, it creates the reservoir and calculates its volume, limited as 11\% of the annual river flow. Next, it compares its size with the water available for storage. Finally, it estimates the costs of the main components of the pumped storage (dam, tunnel, turbine, generator, excavation and land) and its water and energy storage costs for further comparison. The cost estimate is based on charts of a technical manual published by the Norwegian Water Resources and Energy Directorate \cite{SWECO}.

While these references cover a large spectrum of applications using different topologies, none of them translates the identification of hot spots into a mathematical formulation that incorporates the major cost components of PHS. According to Rogeau et al. \cite{ROGEAU2017241}, a cost function would help to eliminate the weakest interconnections and to identify the most economically feasible alternatives regardless of their type of engineering layout. To our knowledge this is the first work with this approach. It uses a DEM to identify reservoirs that meet given storage targets while minimizing embankment, water conveyance and electromechanical equipment (E\&M) costs. Some parameters, such as H and L/H, mentioned above, can be used as criteria to select the most suitable areas to apply this new approach.

\begin{figure}[!ht]
	\centering
	\includegraphics[width=2in]{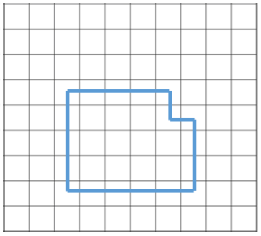}
	\caption{shape example in a grid}
	\label{fig: shape example}
\end{figure}

The aim of this work is to determine the least-cost location and shape of an upper reservoir, as in Figure \ref{fig: shape example} considering a minimum volume requirement and an existing lower reservoir or water body (lake or river), using an integer programming approach. 

The integer programming problem formulation and solution methods are described in section 2. In section 3, the PHS siting and sizing formulation was applied to determine the upper reservoir to be formed on the hillside next to the existing lower reservoir of the Sobradinho hydropower plant. Section 4 presents the modeling results. Finally, section 5 presents the conclusion of this work and future related opportunities.

\section{Proposed formulation}
\subsection{Basic problem}

An integer programming model is proposed for the PHS siting. Binary variables represent grid cells of $n$ line $i$ and $n$ column $j$ that are part of the reservoir perimeter $x_{i,j}$, its interior $y_{i,j}$, or both $z_{i,j}$. Let $\mathcal{X, Y, Z}$ be sets with the allowed perimeter, interior, or both cells. We only allow interior cells that has an elevation bellow the desired head. Perimeter cell must have at least one interior cell as a neighbourhood. A cell is allowed to be a reservoir cell only if it is allowed to be an interior or a perimeter. Cells that already have another reservoir build are removed from this sets. Other reprocessing criteria remove  could also be incorporated into this sets without changing the rest of formulation. Those are all subset of the set $\{1\dots,n\}^2$ where $n$ is the grid length.

Cells allowed to take part of the reservoir are limited by its perimeter, as seen by constraints~\eqref{eq: allowed interior 1}-\eqref{eq: allowed interior 4}. Together, these constraints indicate that if a given cell is part of a reservoir, it must either be part of its boundary (a perimeter cell) or be an interior cell, surrounded by other reservoir cells.

\begin{align}
    z_{i,j} \leq x_{i,j} + z_{i-1,j} \quad \forall i,j \in \mathcal{Z} \label{eq: allowed interior 1} \\
    z_{i,j} \leq x_{i,j} + z_{i+1,j} \quad \forall i,j \in \mathcal{Z} \\
    z_{i,j} \leq x_{i,j} + z_{i,j-1} \quad \forall i,j \in \mathcal{Z} \\
    z_{i,j} \leq x_{i,j} + z_{i,j+1} \quad \forall i,j \in \mathcal{Z} \label{eq: allowed interior 4}
\end{align}

Constraints \eqref{eq: required interior 2} indicate cells that must be in the reservoir given the perimeter. That is, a perimeter cell must have contact with at least three reservoir cell.

\begin{align}
    x_{i,j} \leq z_{i-1,j} + z_{i+1,j} + z_{i,j-1} + z_{i,j+1} \quad \forall i,j \in \mathcal{X} \label{eq: required interior 2}
\end{align} 

Constraints~\eqref{eq: link_variables} impose that a reservoir cell is either in the interior or in the perimeter.

\begin{align}
    z_{i,j} = x_{i,j} + y_{i,j} \quad \forall i,j \in \mathcal{Z} \label{eq: link_variables} \\
\end{align}

Constraints \eqref{eq: interior1}-\eqref{eq: interior4} define interior cells from perimeter and reservoir cells.

\begin{align}
    y_{i,j} \leq z_{i-1,j} \quad \forall i,j \in \mathcal{Y} \label{eq: interior1} \\
    y_{i,j} \leq z_{i+1,j} \quad \forall i,j \in \mathcal{Y} \\
    y_{i,j} \leq z_{i,j-1} \quad \forall i,j \in \mathcal{Y} \\
    y_{i,j} \leq z_{i,j+1} \quad \forall i,j \in \mathcal{Y} \label{eq: interior4}
\end{align}

The minimum volume required is represented by constraint~\eqref{eq: minimum volume}.

\begin{align}
    \sum_{i,j \in \mathcal{Y}} y_{i,j} (H - h_{i, j}) \text{$AreaCell$} \geq \text{$VolMin$} \label{eq: minimum volume}
\end{align}

Where $H$ is a fixed (given) water elevation, $h_{i,j}$ is terrain elevation of cell, \text{$AreaCell$} is the area in each cell, and \text{$VolMin$} is the minimum volume required.


\subsection{Avoiding disconnected reservoirs}

\subsubsection{Traveling salesman problem constraints}

The basic problem formulation can generate multiple disconnected upper reservoirs, as subtours of the traveling salesman problem (TSP) \cite{miller1960integer}. To avoid this, TSP constraints are added using a modified Miller-Tucker-Zemlin (MTZ) formulation are represented by constraints~\eqref{eq: w binary}-\eqref{eq: no_subrouts}. 

\begin{align}
    w_{i,j,h,k} \in \{0,1\} \quad \forall i,j,h,k \in \mathcal{X} \times \mathcal{X}  \label{eq: w binary} \\
    u_{i,j} \in \mathbb{Z}^+ \quad \forall i,j  \in \mathcal{X} \\
    0 \leq u_{i,j} \leq \sum_{h,k} x_{h,k} - 1 \quad \forall i,j  \in \mathcal{X} \\
    0 \leq u_{i,j} \leq x_{i,j} \quad \forall i,j  \in \mathcal{X} \\
    \sum_{h,k} w_{i,j,h,k} = x_{i,j} \quad \forall i,j  \in \mathcal{X} \\
    \sum_{h,k} w_{h,k,i,j} = x_{i,j} \quad \forall i,j  \in \mathcal{X} \\
    u_{i,j} - u_{h,k} + (\sum_{a,b} x_{a,b}) w_{i,j,h,k} \leq \sum_{a,b} x_{a,b} - 1 \quad \forall i,j,h,k \in \mathcal{X} \times \mathcal{X}  \label{eq: no_subrouts}
\end{align}

The modification is that neither the vertex nor the number of vertexes of the graph are an input to the problem. They are decided by variables $x_{i,j}$ and $\sum_{i,j} x_{i,j}$ respectively. Binary variable $w_{i,j,h,k}$ is activated when an edge connects cells $(i,j)$ to $(h,k)$. This can only happen if these cells are neighbors, which requires examining all possible directions. This is not made explicit in the formulation to avoid heavy notation. Therefore, the reader should consider that variables and constraints generations are implicit. Variable $u_{i,j}$ is an auxiliary variable to eliminate subtours.

It is important to remark that constraint~\eqref{eq: no_subrouts} that eliminates subtour in our version of Miller-Tucker-Zemlin formulation contains products of variables. Fortunately, they are products of binary variables, and therefore, linearizable. This linearization is done in the implementation using auxiliary variables.

\subsubsection{Separating planes}

Separating planes can be used to avoid solutions with unconnected reservoirs. If each reservoir were convex, there would no reason for concern. Since that is not the case, they can be used as a tentative to find a quick solution. With a bit of luck, the reservoir will be a connected set. If this is not the case, TSP no-subtour constraints must be added to the model as well. Formulation for the horizontal and vertical planes is presented by Constraints~\eqref{eq: def_variable_up}-\eqref{eq: horizontal_left}.\\

\begin{align}
    up_{i} \in \{0,1\} \quad \forall i  \in \{1,\dots,n\}\label{eq: def_variable_up} \\
    down_{i} \in \{0,1\} \quad \forall i  \in \{1,\dots,n\} \\
    left_{j} \in \{0,1\} \quad \forall j  \in \{1,\dots,n\} \\
    right_{j} \in \{0,1\} \quad \forall j  \in \{1,\dots,n\} 
\end{align}

\begin{align}
    1 - \sum_j y_{i,j} \leq up_{i} + down_{i} \quad \forall i   \in \{1,\dots,n\} \\
    1 - \sum_i y_{i,j} \leq left_{j} + right_{j} \quad \forall j   \in \{1,\dots,n\} 
\end{align}
\begin{align}
    \sum_{k > i, j} y_{k,j} \leq M (1 - down_{i}) \quad \forall i   \in \{1,\dots,n\} \\
    \sum_{k < i, j} y_{k,j} \leq M (1 - up_{i}) \quad \forall i   \in \{1,\dots,n\} \\
    \sum_{i, k > j} y_{i,k} \leq M (1 - left_{j}) \quad \forall j   \in \{1,\dots,n\} \\
    \sum_{i, k < j} y_{i,k} \leq M (1 - right_{j}) \quad \forall j   \in \{1,\dots,n\}  \label{eq: horizontal_left}
\end{align}

\begin{figure}[h!t]
	\centering
	\includegraphics[width=2in]{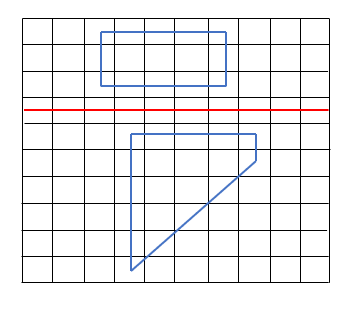}
	\caption{plan view exhibiting a horizontal separating plane}
	\label{fig: Horizontal plane}
\end{figure}

Where $M$ is a number big enough (sometimes referred  as big M) based on the size of the grid. Its function is to allow to turn on/off constraints depending on the value of the auxiliary binary variables. Figure~\ref{fig: Horizontal plane} shows an example of two separate reservoirs and a horizontal separating plane that would avoid this solution.

Diagonals separating planes can be defined as well, similar to horizontal and vertical planes.

\subsection{Objective function}

The three most important cost factors related to PHS facilities in a given location are: (i) construction of embankments; (ii) water conveyance system; (iii) E\&M $c^{e\&m}$. 

Embankments are built to enclose water in reservoirs, avoid seepage whenever the elevation of the terrain in the perimeter of the reservoir is lower than the water elevation. A trapezoidal cross section can be used considering a certain reference width for the crest of the embankment and reference slopes that depend on the embankment material (e.g. a right trapezoid if concrete with a downstream slope of 1V: 0.75H or a isosceles trapezoid with slopes around 1V:2H in case of earth-fill section).  

The cost of water conveyance between the planned upper reservoir and an existing lower reservoir (or river) in the vicinity varies with its length (being a tunnel or a penstock) and the water flow, obtained by $Q=P/(k.H)$, being $k$ the product of water density, acceleration of gravity and equipment efficiency, and thus increasing with the capacity (P) and decreasing with the water head (H). 

The third cost component is related to E\&M. A unit cost can be used $(\$/kW)$ as a function of water head because unit costs tend to decrease for higher heads because the size of turbines, generators, pumps, etc. is reduced while its rotational speed increases.

\subsubsection{Embankment cost}

Assuming a trapezoidal cross-section cost of embankments with a 10m width at the crest and a horizontal : vertical slope of 2:1 for both faces, the cost per cell is given by the cell length of the DEM   ($\text{$Lc$}$) multiplied by the cost of the cross-section, and the unit cost of embankment(\$$5/m^3$), i.e., $c^e_{i,j} = \text{$Lc$}(10(H - h_{i,j}) + 2(H - h_{i,j})^2) $ if $H$ is lower than the terrain elevation of cell $(i, j)$ and zero otherwise.


Note that, in reality, the length and volume, and thus the cost, of the embankment necessary to connect two cells in a diagonal is greater than the values necessary to connect in vertically or horizontally. Auxiliary variables are introduced in the problem to correct this difference in cost. However, they are omitted in the formulation presented here for the sake of simplicity. 

\subsubsection{Water conveyance}

Water conveyance cost is given by water flow $Q = \text{$VolMin$}/\Delta t$. Where \text{$VolMin$} is the desired target volume, and $\Delta t$ is the time it takes to fill up or empty the reservoir. Thus, $Q$ is the desired flow rate.

Additionally, for each possible perimeter cell, one can compute the minimum length $L_{i,j}$ that links it to the other reservoir or river (in the dimensional space) as a pre-processing step prior to optimization. Since, one cannot know which perimeter cell will be connected to the other reservoir or river, a decision link variable is added to the model by adding constraints~\eqref{eq: define_link}-\eqref{eq:  necessary_link}.

\begin{align}
    l_{i,j} \in \{0,1\} \quad \forall i,j \in \mathcal{X} \label{eq: define_link} \\
    l_{i,j} \leq x_{i,j} \quad \forall i,j \in \mathcal{X} \label{eq: limit_link} \\
    \sum_{i,j \in \mathcal{X}} l_{i,j} = 1 \label{eq: necessary_link}
\end{align}

The water conveyance cost is given by two terms: (i) tunnel excavation and treatment; (ii) tunnel lining. The first term is given by a unit cost (e.g. $\$40/m^3)$ x length $L_{i,j}$ ($m$) x cross section $A_{i,j}$ ($m^2$). The later is given by the ratio of the water flow $Q (m^3/s)$ and a maximum admissible velocity (e.g. $4 m/s$). Altogether, the first component can be represented by $c^{ex}_{i,j} = 10QL_{i,j}$.The tunnel lining cost is given by the unit cost of steel (e.g. \$$2.5/kg)$ x total weight ($kg$). Total weight is given by the specific weight of steel ($7840 kg/m^3$) x volume of steel used ($m^3$), which is the tunnel thickness (e.g. 20mm) x perimeter (i.e. $ \pi $ x diameter $D_{i,j}$) x length. If, as before, a maximum admissible velocity is used (e.g. 6 $m/s$) then $Q (m^3/s)$ = 6$A_{i,j}$ = 6$ \pi D_{i,j}^2/4$. Therefore, the diameter is given by $D_{i,j} = \sqrt{2Q/3\pi}$. If steel lining is assumed to be used in 1/3 of the total tunnel length, then the cost is $2.5 \times 7840 \times 0.020 \times \pi \sqrt{2Q/3\pi} \times L_{i,j} / 3$. At last, the lining cost estimate is $c^l_{i,j} =  190 \sqrt{Q} L_{i,j}$

\subsubsection{Equipment costs}

Because the obective of the model is to screen PHS in multiple locations, a general function is used to estimate total E\&M cost. It was derived using a symbolic regression tool based on the Differentiable Genetic Cartesian Programming (dCGP) from the European Space Agency (ESA) \cite{DCGP,DBLP:journals/corr/IzzoBM16} with a dataset constructed with the use of Brazilian Inventory Manual for conventional hydropower projects \cite{MIBRA} for various combinations of water heads, flows and installed capacities. The database is found in the public repository~\cite{andrade_tiago_2020_4019256}.

Worksheets are used in planning and inventory studies to calculate the main dimensions of each hydropower structure and estimate equipment unit costs based on Eletrobras database. They were used to prepare a dataset of total costs for 54 pairs of the following key variables: installed capacity (P, in MW) and hydraulic head (H, in meters). Francis turbines were considered, respecting ranges for P and H. Given that the reference costs were for conventional hydropower, they were  increased by 20\% to consider PHS particularities. 

The machine learning tool dCGP performed a symbolic regression over this dataset, in order to define a simple, yet sensible, component of the objective function that captures the main drivers of total equipment costs and is comparable for multiple locations. The equation selected by the tool was $c^{e\&m} =\dfrac{3068}{\sqrt{H}}+\dfrac{8608}{P}$\ and the final objective function given by \eqref{eq: objective function}.

\begin{align}
    \label{eq: objective function}
    \textit{min} \sum_{i,j \in \mathcal{X} | {H > h_{i,j}}} c^e_{i,j} x_{i,j} + \sum_{i,j} (c^{ex}_{i,j} + c^l_{i,j}) l_{i,j} + {c^{e\&m}P}  
\end{align}

While the first two terms of the objective are intrinsically related the problem formulation, the third component is a constant for a given pair of ${P,H}$. Other cost factors, such as grid connections and access roads can be also considered to compare different alternative sites. Thus, the screening of the best projects consists in solving the optimization problem for different pairs of ${P,H}$ and selecting  alternatives with smallest overall costs.

\subsection{Solution strategies}
\label{solution strategies} 

Due to the computation burden of solving a potentially large integer programming (IP) problem, two strategies are proposed. The first is to simplify the problem by solving first a simplified version of the problem, where only part of the disconnected reservoir constraints are active. If there is only one reservoir, it is best to stop. If not, the alternative is activate the others constraints and try again. The second is to aggregate cells and solve a smaller approximated problem, then, this solution can be used to help solving the larger problem. 

\subsubsection{Disconnected reservoirs strategy}

The different forms of avoiding unconnected reservoirs can be viewed as layers of defense. A possible strategy is to try to solve the problem without any additional constraints. If the solution is not feasible, one can activate the levels of defense in the following recommended order:

\begin{enumerate}
  \item with horizontal and vertical separating planes
  \item with horizontal, vertical, and diagonal separating planes
  \item with no-subtour TSP constraints
\end{enumerate}

\subsubsection{Zoom in heuristic}

The problem size (number of variables and constraints) is proportional to $n^2$, where $n^2$ is the number of cells in a terrain with $n \times n$ cells. Thus, if the DEM has a fine resolution, it may be computationally intensive to solve the original problem and makes sense to simplify the grid by aggregating individual cells into larger ones. The zoom in strategy starts with fewer cells, solves the problem, clips the terrain around the encountered reservoir, increases the precision (reducing cell aggregation) and solves it again. The process is repeated until the finer resolution is achieved (i.e. no aggregation).

\subsection{Other considerations}

There may be cells in the grid that cannot be used due to environmental or any other legal reasons. This can be easily  incorporated in the formulation by logically removing these variables from it. 

A drawback of the formulation is that it contains a product of variables in the TSP constraints. This can be resolved by linearizing the product between binary variables or imposing an upper limit in the number of reservoir perimeter. The latter option can be used to reformulate constraint~\eqref{eq: no_subrouts} by using a big-M instead of the variable products, which appear to be the better option since a natural big-M is the number of squares in the grid. Moreover, the TSP constraints are usually not needed, only the separating planes are enough.

One possible approach to mitigate the issue of adding too many constraints is to use lazy constraints. This could be used both in the separating planes or with the TSP formulation. However, note that the TSP adds variables too. So, even with a lazy constraint implementation, trying to solve the problem first with the separating plans can be advantageous.

Moreover the MTZ TSP formulation has the number of constraints and variables proportional to $n^4$ where $n$ is the length of the grid. We chose this formulation since the problem size is polynomial instead of exponential in the grid length.

To summarize, the model to be solved is to minimize \eqref{eq: objective function} subjective to Constraints \eqref{eq: allowed interior 1}-\eqref{eq: minimum volume},\eqref{eq: def_variable_up}-\eqref{eq: necessary_link} in case the TSP constraints are not necessary. and subject to \eqref{eq: allowed interior 1}-\eqref{eq: necessary_link} otherwise.

\section{Case study}

The PHS siting formulation was applied to determine the upper reservoir to be formed on the hillside next to the existing (lower) reservoir of the Brazilian \textit{Sobradinho} hydropower plant. The selected area, known as \textit{Saco do Arara} had been originally identified as a promising site by former consultant and chairman of the board of The Nature Conservancy, David Harrison \cite{kelman:hal-02147740}, for the following reasons: (i) the proximity with Sobradinho dam decreases costs of connecting the PHS to the power grid; (ii) Sobradinho is a large reservoir with nearly constant water level in a short-term PHS operation; (iii) the elevation of the upper reservoir is at least 150m higher than Sobradinho and the distance between the two reservoirs is not very large. All these factor contribute to reducing water conveyance and E\&M costs; (iv) apparent socio-environmental impacts are also small. Data used for this case study is found in~\cite{andrade_tiago_2020_4019256}.

A field trip to Saco do Arara was made with strong local support by Chesf, the utility that operates Sobradinho and other dams in the São Francisco river of Brazil totalling more than 10 thousand MW in capacity. This support was crucial in terms of logistics given the difficult access to Saco do Arara  site. This field trip confirmed the favorable impression before it was held. 

A DEM of the region with cells with 34m x 34m from NASA's STRM database was considered. It is a square with 266 cells per side ($\approx9km$). From the 70,756 cells of the DEM a total of 19,755 cells (28\%) are part of the Sobradinho reservoir, with an elevation of 385m in the DEM. 

In this study we consider a PHS project with 500 MW of capacity (roughly 50\% of Sobradinho's own capacity) and a storage requirement for 3 hours of operation at this capacity. We also run a sensitivity with a larger stored energy capacity for 12 hours of operation. 

The PHS siting model is executed for three given water heads: 150m, 175m and 200m. In order to identify the effect of the ``zoom in" solution approach, it will be compared with the direct execution of the model (e.g. no heuristic), where variables and constraints of the problem are generated directly for the resolution of the DEM, with no aggregation. 

Because of the combinatorial nature of the problem, the execution may be interrupted - and the best solution so far found exhibited - if the elapsed time reaches the maximum allowed time of 1h. A total of 9 cases are investigated, as shown in Table \ref{tab: tab2}. The computational experiments were done in a i7-8550@1.8GHz with 16 GB, the time limit was 1 hour for all experiments. The models was implemented in Julia \cite{bezanson2017julia} using JuMP \cite{DunningHuchetteLubin2017}, and the IP solver used was Xpress 8.5 \cite{fair2009xpress} with default parameters.

\begin{table}[htbp]
  \footnotesize
  \centering
  \caption{Input data}
    \begin{tabular}{ccrrrrr}
    \toprule
Case  & Head ($m$) & Zoom in & Operation ($h$) & Storage ($hm^3$) \\
              \cmidrule{2-5}
    1     & 150   & yes   & 3     & 5.50 \\
    2     & 175   & yes   & 3     & 4.72 \\
    3     & 200   & yes   & 3     & 4.13 \\
    4     & 150   & no    & 3     & 5.50 \\
    5     & 175   & no    & 3     & 4.72 \\
    6     & 200   & no    & 3     & 4.13 \\
    7     & 150   & yes   & 12    & 22.02 \\
    8     & 175   & yes   & 12    & 18.87 \\
    9     & 200   & yes   & 12    & 16.51 \\
    \bottomrule
    \end{tabular}%
  \label{tab: tab2}%
\end{table}%

They combine those three water heads (150m, 175m and 200m), the use of the \textit{zoom in} heuristic and the required operation time of the PHS: 3h in cases 1-6 (a period associated with peak electricity prices) and 12h in cases 7-9 (a sensitivity to understand the impact on costs of a significant increase in the amount of stored energy). 

The last column of Table \ref{tab: tab2} shows the required storage in $hm^3$. First, it is determined finding the turbined/pumped water flow ($m^3/s$) calculated from the installed capacity and water head, using representative efficiencies, and then multiplying this flow by the duration of the required operation. Notice the reduction of required storage with an increase of head because the turbined/pumped flow is reduced for the same capacity.  

PHS siting cases 1-3 reached the optimum solution. For the remaining cases, the best result found during the allowed computation time is exhibited. Case 6 did not provide a feasible answer within the allowed time of 1 hour. Results of the execution of the PHS siting model are shown in Table \ref{tab: reservoir physical properties}, which is divided into two main groups: key figures for the reservoirs (storage, flooded area and smallest distance to Sobradinho) and constructed embankments (length and volume).

\begin{table}[htbp]
\footnotesize
  \centering
  \caption{PHS siting results: physical properties}
    \begin{tabular}{crrr|rr}
    \toprule
          & \multicolumn{3}{c|}{Reservoir} & \multicolumn{2}{c}{Embankment} \\
   \cmidrule{2-6}
    Case  & Storage ($hm^3$) & Area ($ha$) & \multicolumn{1}{c|}{Distance ($m$)} & Length ($m$) & Volume ($hm^3$) \\
    \cmidrule{2-6}
    1     & 5.56  & 22    & 2729  & 614   & 2.02 \\
    2     & 6.79  & 109   & 1124  & 1380  & 0.26 \\
    3     & 4.13  & 53    & 1206  & 1451  & 0.53 \\
    4     & 5.51  & 13    & 589   & 1973  & 9.51 \\
    5     & 4.73  & 17    & 719   & 2280  & 11.76 \\
    6     & - & -   & -   & -  & - \\
    7     & 22.02 & 86    & 2729  & 2828  & 6.30 \\
    8     & 24.39 & 252   & 886   & 1651  & 0.56 \\
    9     & 16.51 & 137   & 959   & 2383  & 1.55 \\
    \bottomrule
    \end{tabular}%
  \label{tab: reservoir physical properties}%
\end{table}%

Table~\ref{tab: solution_mip} shows the problem size, time elapsed, and achieved gap. All problems, except Case 6 were solved with separating planes. Case 6 with separating reservoirs returned a separable reservoir, thus, we tried to solve it with TSP constraints, but were not able to achieve a feasible solution. To validate the max time, Case 4 was solved again with a 8h time limit, but did not improve the integer solution with the additional time. For the cases which zoom-in strategy was used, reported problem size refers to the larger problem solved for the case and the reported time is the sum of all problems solved. Case 6 size is the problem with TSP constraints.

\begin{table}[htbp]
 \footnotesize
  \centering
  \caption{Optimization summary }
    \begin{tabular}{rrrrr}
    \toprule
    \multicolumn{1}{l}{Case} & \multicolumn{1}{l}{\# Variables} & \multicolumn{1}{l}{\# Constraints} & \multicolumn{1}{l}{time (s)} & \multicolumn{1}{l}{gap (\%)} \\
    \cmidrule{2-5}
     1     & 8168   & 4026   & 19    & 0.0         \\
     2     & 25188  & 10390  & 1654  & 0.0         \\
     3     & 14149  & 6259   & 266   & 0.0         \\
     4     & 160078 & 61672  & {MAX} & 75.3    \\
     5     & 170078 & 65256  & {MAX} & 90.4     \\
     6     & 321657 & 166466 & {MAX} & {-} \\
     7     & 21786  & 9156   & {MAX} & 0.6         \\
     8     & 50360  & 19754  & {MAX} & 7.1         \\
     9     & 28924  & 11758  & {MAX} & 4.1         \\
    \bottomrule
    \end{tabular}%
  \label{tab: solution_mip}%
\end{table}%

\begin{figure}[!ht]
	\centering
	\includegraphics[width=5in]{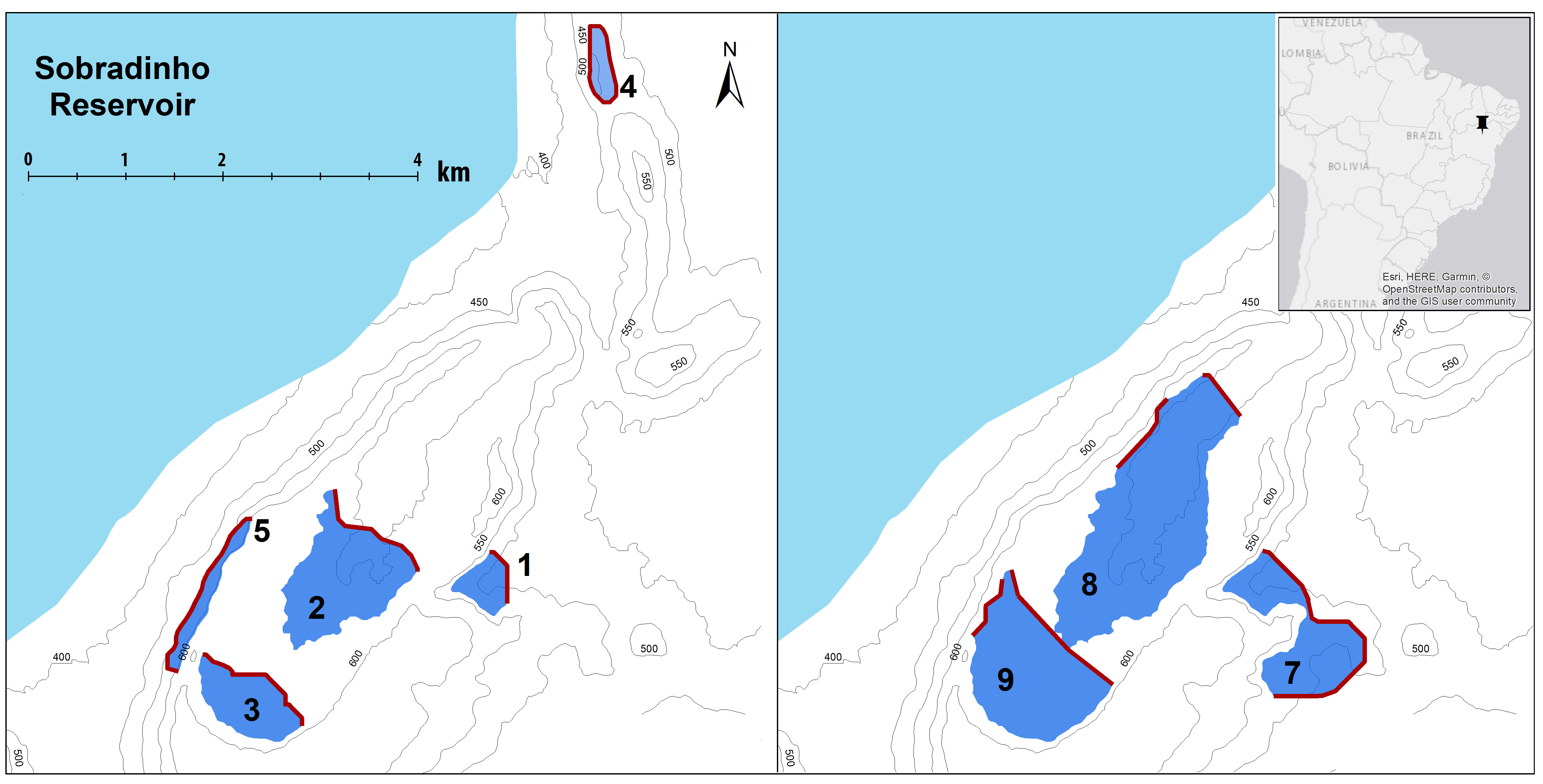}
	\caption{PHS siting optimal reservoirs}
	\label{fig: Cases12345789}
\end{figure}

The first column shows the storage (in million $m^3$) of the upper reservoir. Notice how these values compare with last column of Table \ref{tab: tab2}). In most cases the value is nearly the same, meaning that an increase in storage would increase the project cost. The solution of Case 8, however, has a \textit{surplus} of energy storage of 5.53 $hm^3$ = 24.39 - 18.87. Except for Case 6, the upper reservoirs obtained by the model are shown in the map of Figure~\ref{fig: Cases12345789}.

\section{Discussion}

An examination of Table \ref{tab:costs} shows that costs \textit{decrease} with head not only due to equipment but also, at least in the case study, due to cost savings in the construction of embankments and water conveyance systems. The most important term is the E\&M costs. 

\begin{figure}[!ht]
	\centering
	\includegraphics[width=5.2in]{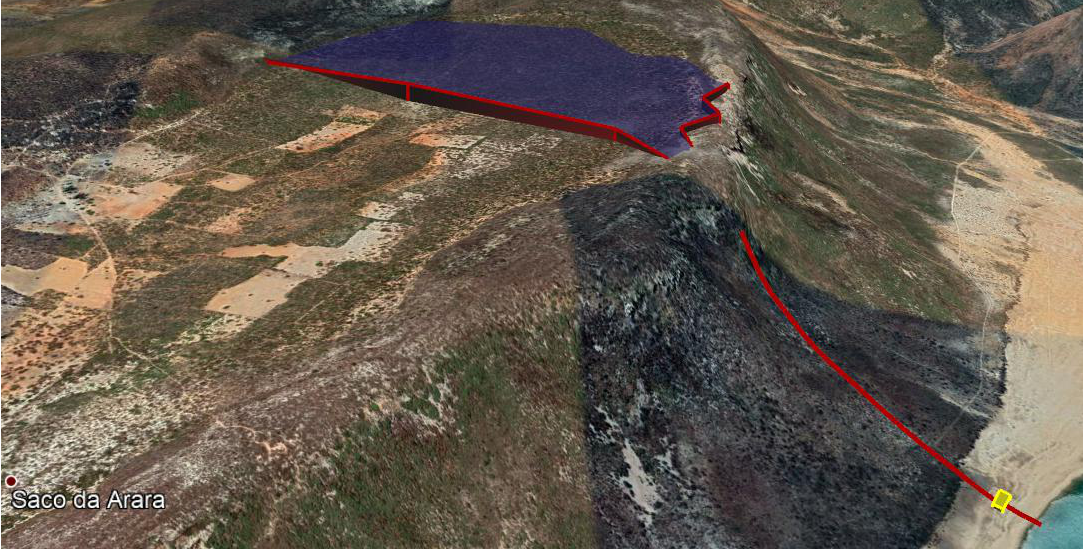}
	\caption{South facing Google Earth scheme of Case 9 with Sobradinho reservoir at lower right near the powerhouse.}
	\label{fig: Google}
\end{figure}

It is interesting to compare the solutions of cases 1-3 with the best feasible solutions (not optimum) found for cases 4-6 because the problems are the same, the difference is just the solution strategy. The comparison of costs is made directly from Table \ref{tab:costs}. Results indicate the effectiveness of the zoom in strategy. Case 2, for instance, was \$55 million cheaper than Case 5. It is also interesting to compare solutions for 3h and 12h of storage. Case 3 and Case 9, whose 3D scheme is plotted in Figure~\ref{fig: Google}, can be used for that matter.

Notice how they use the same approximate location and how the solution of Case 9 contains that of Case 3. The storage difference is 16.51 - 4.13 = 12.38 $hm^3$ and the cost increase \$131 - \$128 = \$3 million. Thus, the \textit{incremental} cost of storage is 3/12.38 $\approx$ $0.24/m^3$, which translates into \$0.5/kWh of storage. This value is two to three orders of magnitude cheaper than the present costs of Lithium-Ion batteries. 

This low incremental storage costs is typical of PHS projects. This raises the following question: how much storage should one select? The actual target must come from top-down studies \cite{kelman:hal-02147740} that investigate the role of PHS in the power systems to provide support for variable renewable energy sources, such as wind power and solar photovoltaic, to reduce the need to expand transmission lines for energy imports in some areas of the grid, or others grid-dependent economic benefits. 

\begin{table}[htbp]
 \footnotesize
  \centering
  \caption{PHS siting: Costs (\$ Millions) }
    \begin{tabular}{rrrrr}
    \toprule
    \multicolumn{1}{l}{Case} & \multicolumn{1}{l}{Embankment} & \multicolumn{1}{l}{Water Conveyance} & \multicolumn{1}{l}{E\&M Equipment} & \multicolumn{1}{l}{Total} \\
    \cmidrule{2-5}
    1     & 10    & 22    & 134   & 166 \\
    2     & 1     & 8     & 125   & 134 \\
    3     & 3     & 8     & 117   & 128 \\
    4     & 48    & 5     & 134   & 187 \\
    5     & 59    & 5     & 125   & 189 \\
    6     & -    & -     & -   & - \\
    7     & 31    & 22    & 134   & 187 \\
    8     & 3     & 6     & 125   & 134 \\
    9     & 8     & 6     & 117   & 131 \\
    \bottomrule
    \end{tabular}%
  \label{tab:costs}%
\end{table}%

\section{Conclusions and future works}

Given the seasonal and variable electricity production of renewable sources, and the interest of maintaining a balance between supply and demand, energy storage strategies can be included as a resource in integrated planning studies of power systems. By far, the largest technology used globally to this end is PHS because of the fast response of power, relatively large storage capacity and economic competitiveness. 

The contribution of this paper is the formulation of a PHS siting model composed of two main parts. In the first one, it was shown how constraints determine the reservoir which must meet a given target storage. In the second part, different solution strategies were discussed to eliminate unconnected reservoirs from the solution. 

Further research can be developed to facilitate a step of this PHS screening process, prior to the optimization model developed here. This first screening should identify bow-shaped locations in larger areas to reduce the terrain area that is inputted in the model. Siting model using modern computational techniques in GIS-based environments as some references mentioned in this paper will work well to this end.

A case study for the \textit{Saco do Arara} site was conducted by changing key parameters in the PHS siting model, and maintaining the digital elevation model. This can be interpreted as a \textit{bottom-up} approach where to define competitive PHS projects. The model is run several times for different locations based on local DEMs that were, themselves, the result of the application filters on much larger areas as described in references to find \textit{hot spots}. 

Optimum solutions are found and the best results, when various locations are considered, must together meet grid requirements coming from the \textit{top-down} approach based, for instance, on the results of an integrated resource planning model.

A possible future development is integrate this formulation with a \textit{top-down} approach that determines targets for energy storage based on an integrated planning of the power systems. Another possibility is to integrate the model with a more elaborate engineering design module to improve the cost estimates. This approach is similar to the approach used in the Brazilian Inventory Manual. Finally the model could be deployed in GIS-based platforms, as in the case of STORES, to facilitate the visualization and interpretation of the model results.

Finally, while it is possible to expand the model to incorporate other cost terms, it is quite difficult to stipulate the costs of socio-environmental impacts. One approach would include: (i) the use of general layers of information, such as protected areas, priority areas for conservation, land cover and traditional communities; (ii) these layers should be used to generate a measure of the impacts of each candidate PHS project through metrics; (iii) these metrics would be associated with additional costs, whether direct (e.g. land acquisition) or indirect (e.g. compensations); (iv) the formulation would add these cost terms to the objective function being minimized. Procedure (i)-(iv) would differentiate projects as part of the screening process.

\bibliography{myref.bib}

\bibliographystyle{IEEEtran}

\end{document}